\documentclass[12pt]{article}
\usepackage{amsmath,amsfonts,amssymb,amsthm} 
\newtheorem{theorem}{Theorem}[section]
\newtheorem{prop}[theorem]{Lemma}
\newtheorem{coro}[theorem]{Corollary}
\newtheorem{prob}{Problem}
\newtheorem*{laglem}{Lagrange Lemma}
\newtheorem*{conj}{Conjecture}
\theoremstyle{remark}
\newtheorem{remark}[theorem]{Remark}
\numberwithin{equation}{section}
\newcommand\la{\lambda}

\newcommand\nbi[3]{{\binom{#1}{#2}}_{#3}}
\newcommand\pbi[2]{\genfrac{\langle}{\rangle}{0pt}{}{#1}{#2}}
\newcommand\npbi[3]{{\genfrac{\langle}{\rangle}{0pt}{}{#1}{#2}}_{#3}}

\title{Jack polynomials and some identities for partitions}
\author{Michel Lassalle\\
\small Centre National de la Recherche Scientifique\\[-0.8ex]
\small Institut Gaspard Monge, Universit\'e de Marne-la-Vall\'ee\\[-0.8ex]
\small 77454 Marne-la-Vall\'ee Cedex, France\\[-0.8ex]
\small \texttt{lassalle @ univ-mlv.fr}\\[-0.8ex]
\small \texttt{http://igm.univ-mlv.fr/{\textasciitilde}lassalle/index.html}}
\date{\small {\em Mathematics Subject Classification 2000 :} 
05A10, 05A17, 05E05, 33C52, 33C80}

\begin{document}
\maketitle
\begin{abstract}
We prove an identity about partitions involving new combinatorial 
coefficients. The proof given is using a generating function. As an application 
we obtain the explicit expression of two shifted symmetric functions, 
related with Jack polynomials. These quantities 
are the moments of the ``$\alpha$-content'' random variable 
with respect to some transition probability distributions.
\end{abstract}

\section{Introduction}

A partition $\la= (\la_1,...,\la_n)$ 
is a finite weakly decreasing
sequence of positive integers, called parts. The number 
$n=l(\la)$ of parts is called the length of 
$\la$, and
$|\la| = \sum_{i = 1}^{n} \la_i$
the weight of $\la$. For any $i\geq1$, 
$m_i(\la)  = \textrm{card} \{j: \la_j  = i\}$
is the multiplicity of $i$ in $\la$.  Clearly one has
$l(\la)=\sum_{i} m_i(\la)$ and $|\la|=\sum_{i} i m_i (\la)$. 
We identify $\la$ with its Ferrers diagram
$\{ (i,j) : 1 \le i \le l(\la), 1 \le j \le \la_i \}$ 
and set 
\[z_\la  = \prod_{i\ge1} 
{i}^{m_i(\la)} m_i(\la) ! .\]

In this paper we consider the problem of evaluating 
sums of the following type
\[\sum_{|\mu| = n} \frac{1}{z_\mu}\prod_{k \ge 1}
 (S_k)^{m_k(\mu)},\]
where $S_k$ is a formal series, depending on some indeterminates.
The more elementary sum of this kind is well known~\cite{Ma}. 
Let $X$ be an indeterminate and $n$ a positive integer. We write 
\[(X)_n = X(X+1) \ldots (X+n-1) \quad,\quad 
[X]_n = X(X-1) \ldots (X-n+1)\]
for the raising and lowering factorials and
$\binom{X}{n}  = [X]_n/n!$. Choosing $S_k=X$, one has
\[\sum_{|\mu| = n} \frac{X^{l(\mu)}}{z_\mu}  = 
\binom{X + n - 1}{n} = \frac{(X)_{n}}{n!}.\]	

Let $X_{0}$ and $X=(X_{1},X_{2},\ldots)$ be 
(infinitely many) independent indeterminates. 
In~\cite{LL} we considered another summation 
of the same type
\[\sum_{|\mu| = n} \frac{1}{z_\mu}\prod_{k \ge 1}
 {\left(\sum_{r \ge 0}u^r\frac{(k)_r}{r!} X_r
 \right)}^{m_k(\mu)}\]
and gave its formal series expansion.
 
In this paper we investigate a natural generalization of 
our previous result and give the formal series expansion of
\[\sum_{|\mu| = n} \frac{1}{z_{\mu}}\prod_{k \ge 1}
{\left(\sum_{r,s \ge o} u^r v^s \frac{(k)_r}{r!}\,\frac{(k)_s}{s!}\,
 X_{r+s}\right)}^{m_{k}(\mu)}.\] 

The interest of this result is twofold. Firstly it involves new combinatorial objects 
whose background remains quite mysterious. Actually the given explicit expansion 
is written in terms of a new family 
of coefficients associated to partitions, which are themselves built 
in terms of another new family of positive integers. Both objects 
generalize classical binomial coefficients, but it
should be emphasized that their combinatorial 
interpretation remains very obscure. It is an intriguing problem to clarify 
this underlying structure, which would at the same time
provide a bijective proof of our result.

Secondly any specialization of the indeterminates 
$X_i$ may lead to a different application of our formula. 
Here our attention keeps mainly focused on symmetric function theory and 
(shifted) Jack polynomials.

Let $\alpha$ be some positive real number, $\la$ a partition and $A_\la$ 
the family of its so called ``$\alpha$-contents'' $\{j-1-(i-1)/\alpha, 
(i,j) \in \la\}$. Specializing  $X_i$ to the $i$th power sum symmetric 
function $p_i$ evaluated on $A_\la$, we obtain the series expansion of
\[\frac{(x+y+1)_{\la}}{(x+y)_{\la}}\, \frac{(x)_{\la}}{(x+1)_{\la}},\]
where $x,y$ are two indeterminates, and $(x)_\la$ is some natural generalization 
of the classical raising factorial.

It turns out that two families of shifted 
symmetric functions, closely related to (shifted) Jack 
polynomials, have a generating function of this type. 
Our formula provides an explicit expression for these 
shifted symmetric functions, thus proving (and improving) some earlier 
conjectures~\cite{La4,La5}. These expressions are new even 
in the case $\alpha=1$, 
which corresponds to (shifted) Schur functions. 

Our results have some connection with the works 
of Kerov~\cite{Ke,Ke2}. To any
Young diagram are associated two discrete probability 
distributions, generalizing the Plancherel 
transition and co-transition classical probabilities. It turns out that 
the two shifted symmetric functions explicited here are the 
moments of the $\alpha$-content random variable with respect
to these distributions.

Another consequence of our results is to allow a very easy 
computation of the expansion of Jack polynomials in 
terms of power sum symmetric functions. We hope to report 
about this application in a forthcoming paper.
 
\section{Combinatorial tools}

\subsection{Positive integers}

We first recall some results about a new family of positive integers, 
which we have introduced in~\cite{La1}. 

Let $n,p,k$ be three integers with $1 \le k \le n$ 
and $0 \le p \le n$. We define
\begin{equation*}
{\binom{n}{p}}_{k} 
= \frac{n}{k} \ \sum_{r \ge 0} 
\binom{p}{r}\binom{n-p}{r} 
\binom{n-r-1}{k-r-1}.
\end{equation*}
We have obviously
\[{\binom{n}{p}}_{k} = 0 \quad \textrm{for} 
\quad k>n \quad , \quad 
\nbi{n}{p}{1}=n  
\quad , \quad  
{\binom{n}{p}}_{k} ={\binom{n}{n-p}}_{k} .\]
These numbers generalize 
the classical binomial coefficients, since we have
\[\nbi{n}{0}{k} = \binom{n}{k} \quad,\quad \nbi{n}{1}{k}=k \binom{n}{k}
\quad,\quad \nbi{n}{p}{n} = \binom{n}{p},\]
the last property being a direct consequence of the classical 
Chu-Vandermonde formula.

Other special values have been computed in~\cite{La1}, where 
it was also shown that the numbers $\nbi{n}{p}{k}$ are 
\textit{positive integers}. A combinatorial interpretation of
these numbers has been given in \cite{JLZ}.

In~\cite{La1} we gave the following generating function
\begin{equation}
\sum_{k= 1}^n {\binom{n}{p}}_{k} 
{\left(\frac{z}{1-z}\right)}^k =
 n z \;{}_2F_1\!\left[
\begin{matrix}p+1,n-p+1\\
2\end{matrix};z\right],
\end{equation}
where as usual the classical Gauss hypergeometric function is denoted by
\[{}_2F_1\!\left[\begin{matrix}a,b\\
c\end{matrix};z\right]=
\sum _{i \ge 0}\frac {(a)_i(b)_i}
{(c)_i}\;\frac{z^i}{i!}.\] 
In other words the generating function
\[G_n(y,z) =\sum_{p=0}^{n} \sum_{k=1}^{n} \nbi{n}{p}{k} 
y^p {\left(\frac{z}{1-z}\right)}^k\]
can be written
\begin{equation}
G_n(y,z) =
n z \sum_{p=0}^n y^p  \,{}_2F_1\!\left[\begin{matrix}p+1,n-p+1\\
2\end{matrix};z\right].
\end{equation}
It must be emphasized that $\nbi{n}{p}{k}$ is 
\textit{not} defined for $p>n$, and that the generating 
function (2.1) is only valid for $0 \le p \le n$. 

An explicit expression is known for 
$G_n(y,z)$~\cite{Z,La1}, but we shall not need it. With $x=z/(1-z)$, 
i.e. $z= x/(1+x)$, we have
\begin{multline*}
G_n(y,z)= 2^{-n} {\big( (1+x)(1+y)+\sqrt{(1+x)^2(1+y)^2-4y(1+x)} 
\big)}^n\\
+ 2^{-n} {\big( (1+x)(1+y)-\sqrt{(1+x)^2(1+y)^2-4y(1+x)} \big)}^n -1-y^n.
\end{multline*}

\subsection{New coefficients}

The following notion has been introduced in previous 
papers~\cite{La2,La3}. For $\la$ a partition and any 
integer $k \ge 1$, 
let $\pbi{\la}{k}$ denote the number of ways to choose $k$ different cells 
in the Ferrers diagram of $\la$, taking \textit{at least one cell from each row}.
It is easily seen that
\[\pbi{\la}{k}= \sum_{(k_i)}\prod_{i=1}^{l(\la)} \binom{\la_i}{k_i}\]
the sum being taken over all decompositions $k=\sum_{i=1}^{l(\la)} k_i$ with $k_i 
\neq 0$ for any $i$. The generating function for $\pbi{\la}{k}$ is the following
\begin{equation*}
\begin{split}
\sum_{r \ge  1} \pbi{\la}{r}\,x^{r} &= 
\prod_{i=1}^{l(\la)} \left({(1 + x)}^{\la_i} - 1\right)\\
 &= \prod_{i \ge  1} {\left({{(1 + x)}^i - 1}\right)}^{m_i(\la)}.
\end{split}
\end{equation*}

Now in a strictly parallel way, for any integers $0 \le p \le |\la|$ and 
$k \ge 1$, we define
\[\npbi{\la}{p}{k}= \sum_{(p_i)} \sum_{(k_i)} \prod_{i=1}^{l(\la)} 
\nbi{\la_i}{p_i}{k_i}\]
the sum being taken over all decompositions $p=\sum_{i=1}^{l(\la)} p_i$, 
$k=\sum_{i=1}^{l(\la)} k_i$ with $0 \le p_i \le \la_i$ and 
$k_i \neq 0$ for any $i$. 
Observe that there is no such restriction for $p_i$. 

This definition yields easily
\[\npbi{\la}{p}{k}=0 \quad \textrm{except if} \quad l(\la) \le k \le |\la|.\]
Indeed it is obvious that $\npbi{\la}{p}{k}=0$ for $k<l(\la)$, and since 
$\nbi{n}{p}{k} = 0$ for $k>n$, we have also $\npbi{\la}{p}{k}=0$ for
$k>|\la|$.
For instance $\npbi{\la}{p}{1}=0$ except if $\la$ is a row partition $(n)$. 
In this case we have $\npbi{(n)}{p}{k}=\nbi{n}{p}{k}$.

We have obviously
\[\npbi{\la}{0}{k}= \pbi{\la}{k} \quad, \quad \npbi{\la}{1}{k}= k 
\pbi{\la}{k},\]
and also
\[\npbi{\la}{p}{k}=\npbi{\la}{|\la|-p}{k}.\]
Finally it is not difficult to check that
\[\npbi{\la}{p}{|\la|}=\binom{|\la|}{p}
\quad,\quad \npbi{\la}{p}{|\la|-1}=(|\la|-m_1(\la))
\left[\binom{|\la|-1}{p}+\binom{|\la|-2}{p-2}\right].\]
For instance the first relation is a direct consequence of the 
Chu-Vandermonde formula. Indeed the definition easily implies 
\[\npbi{\la}{p}{|\la|}= \sum_{(p_i)}\prod_{i=1}^{l(\la)} 
\nbi{\la_i}{p_i}{\la_i} = \sum_{(p_i)}\prod_{i=1}^{l(\la)} 
\binom{\la_i}{p_i}.\]

As a straightforward consequence of their definition, 
the generating function for the numbers $\npbi{\la}{p}{k}$ is the following
\begin{equation}
\sum_{p=0}^{|\la|} \sum_{k=l(\la)}^{|\la|} 
\npbi{\la}{p}{k} y^p {\left(\frac{z}{1-z}\right)}^k 
=\prod_{i=1}^{l(\la)} G_{\la_i}(y,z)
=\prod_{i \ge 1} {\Big(G_i(y,z) \Big)}^{m_i(\la)}.
\end{equation}
One has easily $G_i(0,\frac{x}{1+x})=(1+x)^i-1$. Thus for $y=0$ we recover the generating 
function of $\npbi{\la}{0}{k}= \pbi{\la}{k}$.

It must be emphasized that $\npbi{\la}{p}{k}$ is 
\textit{not} defined for $p>|\la|$. 
An interesting problem is to get a combinatorial interpretation 
of $\npbi{\la}{p}{k}$.

\subsection{Auxiliary polynomials}

Let $X=(X_{1},X_{2},\ldots)$ be (infinitely many) independent indeterminates. 
In~\cite{LL,La1,La2} for any integers $n,k \ge 1$ we have defined
\[P_{nk}(X)=\sum_{|\mu| = n} 
\frac{\pbi{\mu}{k}}{z_{\mu}}\prod_{i \ge 1}{X_{i}}^{m_{i}(\mu)}.\]
For $k=0$ we got $P_{n0}(X)=0$ ($n \neq 0$), and $P_{00}(X)=1$.

In a strictly parallel way, for any integers 
$n \ge 1$, $k \ge 1$ and $0 \le p \le n$ we set
\[P_{npk}(X)=\sum_{|\mu| = n} 
\frac{\npbi{\mu}{p}{k}}{z_{\mu}}\prod_{i \ge 1}{X_{i}}^{m_{i}(\mu)}.\]
Since $\npbi{\mu}{p}{k}=0$ for $k<l(\mu)$, this sum is restricted to partitions 
such that $l(\mu) \le k$. Hence $P_{npk}(X)$ is a polynomial in $X$. 
Similarly since $\npbi{\mu}{p}{k}=0$ for $k>|\mu|$,
one has $P_{npk}(X)=0$ for $k>n$. 

For $k=0$ it is natural to extend the previous 
definition by the convention
$P_{np0}(X)=0$ with the only exception $P_{000}(X)=1$.
It must be emphasized that $P_{npk}(X)$ is \textit{not} 
defined for $p > n$.
 
For $k=1$ we have $P_{np1}(X)=X_n$ and for $k=n$
\[P_{npn}(X)= \binom{n}{p} \sum_{|\mu| = n}
\frac{1}{z_{\mu}}\prod_{i \ge 1}{X_{i}}^{m_{i}(\mu)}=\binom{n}{p}P_{nn}(X).\]
Finally we have $P_{n0k}(X)= P_{nk}(X)$ and $P_{n1k}(X)= 
kP_{nk}(X)$. It is also obvious that 
$P_{npk}(X)=P_{n,n-p,k}(X)$.

\section{Main identity}

Let $X_{0}$ and $X=(X_{1},X_{2},\ldots)$ be (infinitely many) 
independent indeterminates. We give the explicit evaluation of 
\[\sum_{|\mu| = n} \frac{1}{z_\mu}\prod_{k \ge 1}
 (S_k)^{m_k(\mu)},\]
when the formal series $S_k$ is chosen to be 
\[S_k=\sum_{r,s \ge 0} u^r v^s \frac{(k)_r}{r!}\,\frac{(k)_s}{s!}\,
  X_{r+s}.\]

\begin{theorem}
For any integer $n \ge 1$ we have
\begin{multline*}
\sum_{|\mu| = n} \frac{1}{z_{\mu}}\prod_{k \ge 1}
 {\left(\sum_{r,s \ge 0} u^r v^s \frac{(k)_r}{r!}\,\frac{(k)_s}{s!}\,
  X_{r+s}\right)}^{m_{k}(\mu)} =\\
\sum_{p,q \ge 0}  u^p v^q
\left( \sum_{k=0}^{\mathrm{min}(n,p+q)} \binom{X_{0}+n-1}{n-k} P_{p+q,p,k}(X) 
\right).
\end{multline*}
\end{theorem}

\begin{proof}
Denote by $L(n)$ (resp. $R(n)$) the left (resp. right)-hand side of this 
identity. On the right-hand side, with $b=a/(1-a)$ and $w=u/v$, we get
\begin{equation*}
\begin{split}
\sum_{n \ge 1} a^n R(n) &=
\sum_{p,q \ge 0} u^p v^q  \sum_{k = 0}^{p+q} \left( \sum_{n \ge k} 
a^n \binom{X_0+n-1}{n-k} \right) P_{p+q,p,k}(X)  \\ 
& = (1-a)^{-X_0} \sum_{p,q \ge 0} u^p v^q \sum_{k = 0}^{p+q} b^k
\left( \sum_{|\mu|=p+q} 
\frac{\npbi{\mu}{p}{k}}{z_{\mu}}\prod_{i \ge 1}{X_i}^{m_i(\mu)} \right)\\
& = (1-a)^{-X_0} \sum_{\mu} v^{|\mu|} \sum_{p = 0}^{|\mu|} \sum_{k = 0}^{|\mu|} 
w^p b^k \frac{\npbi{\mu}{p}{k}}{z_{\mu}}\prod_{i \ge 1}{X_i}^{m_i(\mu)}.
\end{split}
\end{equation*}
Using the generating function (2.3) for $\npbi{\mu}{p}{k}$, we obtain
\begin{equation*}
\begin{split}
\sum_{n \ge 1} a^n R(n) &=
(1-a)^{-X_0} \sum_{m_1,m_2,\ldots \ge 0}
\prod_{k \ge 1} \frac{1}{m_k!} {\left( v^k \frac{X_k}{k} \ G_k(w,a) 
\right)}^{m_k}\\
&=(1-a)^{-X_0} \prod_{k \ge 1} \textrm{exp} \left( v^k \frac{X_k}{k} \
G_k(w,a) \right)\\
&=(1-a)^{-X_0}  \textrm{exp} \left( \sum_{k \ge 1} v^k 
\frac{X_k}{k} \ G_k(w,a) \right).
\end{split}
\end{equation*}
On the left-hand side of the identity we have 
\begin{equation*}
\begin{split}
\sum_{n \ge 1} a^n L(n) &=
\sum_{\mu} \frac{a^{|\mu|}}{z_{\mu}}\prod_{k \ge 1}
 {\left(X_{0}+\sum_{t \ge 1} X_{t} \,
 \sum_{r=0}^t u^r v^{t-r} \frac{(k)_r}{r!} \, \frac{(k)_{t-r}}{(t-r)!}  
 \right)}^{m_{k}(\mu)}\\
 &= \sum_{m_1,m_2,\ldots \ge 0}
\prod_{k \ge 1} \frac{1}{m_k!} 
{\left[\frac{a^k}{k} \left(X_{0}+\sum_{t \ge 1} v^t X_{t} \,
 \sum_{r=0}^t w^r \frac{(k)_r}{r!}\, \frac{(k)_{t-r}}{(t-r)!} 
 \right)\right]}^{m_k}\\
&=(1-a)^{-X_0} \prod_{k \ge 1} \textrm{exp} \left(\frac{a^k}{k} 
\sum_{t \ge 1} v^t X_{t} \,
 \sum_{r=0}^t  w^r \frac{(k)_r}{r!}\, \frac{(k)_{t-r}}{(t-r)!} \right)\\
 &=(1-a)^{-X_0} \textrm{exp}\left[ \sum_{t \ge 1} v^t X_{t} \, 
 \sum_{r=0}^t w^r \left( \sum_{k \ge 1} \frac{a^k}{k} \, \frac{(k)_r}{r!} 
 \, \frac{(k)_{t-r}}{(t-r)!} \right) \right].
\end{split}
\end{equation*}
Thus it is enough to prove
\begin{equation*}
G_n(w,a) = n \sum_{r=0}^n w^r \left( \sum_{i \ge 1} 
\frac{a^i}{i}\, \frac{(i)_r}{r!} \, \frac{(i)_{n-r}}{(n-r)!} \right).
\end{equation*}
But obviously
\[\sum_{i \ge 1} \frac{a^i}{i}\, \frac{(i)_r}{r!} \,
 \frac{(i)_{n-r}}{(n-r)!} = a \,{}_2F_1\!\left[\begin{matrix}r+1,n-r+1\\
2\end{matrix};a\right].\]
We conclude by applying relation (2.2).
\end{proof}

Theorem 3.1 immediatly yields
\begin{multline*}
\sum_{|\mu| = n} (-1)^{n-l(\mu)}\frac{1}{z_{\mu}}\prod_{k \ge 1}
 {\left(\sum_{r,s \ge 0} u^r v^s \frac{(k)_r}{r!}\,\frac{(k)_s}{s!}\,
  X_{r+s}\right)}^{m_{k}(\mu)} =\\
\sum_{p,q \ge 0}  u^p v^q
\left( \sum_{k=0}^{\mathrm{min}(n,p+q)} 
{(-1)}^k\binom{X_{0}-k}{n-k} P_{p+q,p,k}(-X) 
\right).
\end{multline*}
It is an open question whether the right-hand side may be expressed 
in terms of the quantities 
$P_{p+q,p,k}(X)$. This is the case for $v=0$ (\cite{LL}, Theorem 
2, p. 301). Then one has
\begin{equation*}
\sum_{|\mu| = n} (-1)^{n-l(\mu)} \frac{1}{z_\mu} \prod_{k \ge 1}
 {\left(\sum_{r \ge 0}u^r\frac{(k)_r}{r!} X_r
 \right)}^{m_k(\mu)} =
\sum_{p \ge 0} u^p\left( \sum_{k=0}^{\mathrm{min}(n,p)} \binom{X_0-p}{n-k} 
P_{pk}(X) \right).
\end{equation*}
This is a consequence of the following non 
trivial property (\cite{JZ}, Corollaire 5),
\[\sum_{k=l(\mu)}^{\mathrm{min}(n,|\mu|)} \binom{X_0-|\mu|}{n-k} 
\pbi{\mu}{k}=
\sum_{k=l(\mu)}^{\mathrm{min}(n,|\mu|)} (-1)^{k-l(\mu)} 
\binom{X_0-k}{n-k} \pbi{\mu}{k}.\]
It would be interesting to obtain a generalization of this 
identity for $\npbi{\mu}{p}{k}$.

\section{First specialization}

Various applications of Theorem 3.1 may be obtained by using different
specia-lizations of the indeterminates $X_i$. In this section we shall 
consider the simplest case, obtained when all the $X_i, i\ge 0$ are equal.

Let $x$ be an indeterminate. In~\cite{La3} (Theorem 1, 
p. 461) we gave a bijective proof 
(due to Rodica Simion) of the identity
\[\sum_{|\mu| = n} \pbi{\mu}{k}\frac{x^{l(\mu)}}{z_\mu}=
\binom{n-1}{k-1}\binom{x + k - 1}{k}.\]	
Since the Stirling numbers of the first kind $s(k,r)$ are defined by 
the generating function
\[(x)_k =\sum_{r\ge1} |s(k,r)| \,x^r,\]
this identity is equivalent to
\[\binom{n-1}{k-1} |s(k,r)|= k! \, 
\sum_{\begin{subarray}{c}|\mu|=n\\l(\mu) = r\end{subarray}} 
\frac{\pbi{\mu}{k}}{z_{\mu}}.\]

We generalize this result as follows (the previous case
may be recovered by choosing $p=0$).
\begin{theorem} Let $x$ be an indeterminate. For any nonnegative integers 
$n,p,r$, we have
\[\sum_{|\mu| = n} \npbi{\mu}{p}{k}\frac{x^{l(\mu)}}{z_\mu}=
\frac{k}{n}\nbi{n}{p}{k}\binom{x + k - 1}{k}.\]
\end{theorem}
\begin{proof}
With the specialization $X_i=x$ for any $i\ge 0$, the left-hand side 
of the identity of Theorem 3.1 reads
\begin{equation*}
\begin{split}
\sum_{|\mu| = n} \frac{1}{z_{\mu}}\prod_{k \ge 1}
 {\left(\sum_{r,s \ge 0} u^r v^s \frac{(k)_r}{r!}\,\frac{(k)_s}{s!}\, 
 x \right)}^{m_{k}(\mu)} &=
\sum_{|\mu| = n} \frac{x^{l(\mu)}}{z_{\mu}}
\prod_{k \ge 1}{\left((1-u)^{-k}\,(1-v)^{-k}\, \right)}^{m_{k}(\mu)}\\
&=\binom{x+n-1}{n} (1-u)^{-n}\,(1-v)^{-n}\\
&=\binom{x+n-1}{n} \left(\sum_{p,q \ge 0} u^p v^q 
\frac{(n)_p}{p!}\,\frac{(n)_q}{q!} \right).
\end{split}
\end{equation*}
By identification of the coefficients in $u$ and $v$, we immediatly 
obtain for $p,q > 0$,
\begin{equation*}
\begin{split}
\binom{x+n-1}{n} \frac{(n)_p}{p!}\,\frac{(n)_q}{q!} &=
\sum_{k=1}^{\mathrm{min}(n,p+q)} \binom{x+n-1}{n-k} P_{p+q,p,k}(x,\ldots,x) \\
&=\sum_{k=1}^{\mathrm{min}(n,p+q)} \binom{x+n-1}{n-k}
\left(\sum_{|\mu| = p+q} \npbi{\mu}{p}{k}\frac{x^{l(\mu)}}{z_\mu}\right).
\end{split}
\end{equation*}
On the other hand, it is easily proved that we have
\begin{equation*}
\binom{x+n-1}{n} \frac{(n)_p}{p!}\,\frac{(n)_q}{q!} 
=\sum_{k=1}^{\mathrm{min}(n,p+q)} \binom{x+n-1}{n-k}
\frac{k}{p+q}\nbi{p+q}{p}{k}\binom{x + k - 1}{k},
\end{equation*}
equivalently
\[\frac{(n)_p}{p!}\,\frac{(n)_q}{q!} 
=\sum_{k=1}^{\mathrm{min}(n,p+q)} \binom{n}{k} \frac{k}{p+q} \nbi{p+q}{p}{k}.\]
Indeed for $p=0$ this is the classical Chu-Vandermonde formula.
The property is then a direct consequence of 
the following recurrence relation (\cite{La1}, Lemma 3.2)
\[p{\binom{p+q}{p}}_k
=(q+1){\binom{p+q}{p-1}}_k -\frac{p+q}{p+q-1}(q-p+1){\binom{p+q-1}{p-1}}_k.\]
But now if for a sequence of functions $f_k^{(pq)}(x),1 \le k\le p+q$, the relation
\[\sum_{k=1}^{\mathrm{min}(n,p+q)} \binom{x+n-1}{n-k} f_k^{(pq)}(x)=0\]
is satisfied for any $n \ge 1$, by choosing successively 
$n=1,2,\ldots,p+q$ one has iteratively $f_k^{(pq)}=0$, and we can conclude.
\end{proof}

Theorem 4.1 can be equivalently stated as
\[\frac{k}{n}\nbi{n}{p}{k} |s(k,r)| = k! \, 
\sum_{\begin{subarray}{c}|\mu|=n\\l(\mu) = r\end{subarray}} 
\frac{\npbi{\mu}{p}{k}}{z_{\mu}}.\]
It would be interesting to obtain a bijective proof, in the same 
vein than Simion's.

\section{Main specialization}

We now consider another specialization of Theorem 3.1 
in the framework of (shifted) Jack polynomials.

Let $\alpha$ be some positive real number, and $\la$ an 
arbitrary partition. For any cell $(i,j)$ of $\la$ we define 
the ``$\alpha$-content'' 
\[c_{ij}^{\alpha}= j-1-(i-1)/\alpha.\]
We consider the following natural generalization of the ``raising'' 
and ``lowering factorial''
\[(x)_{\la}=\prod_{(i,j) \in \la} \left(x+c_{ij}^{\alpha}\right)
\quad,\quad
[x]_{\la}=\prod_{(i,j) \in \la} \left(x-c_{ij}^{\alpha}\right).\]
For any integer $k\ge 0$ we define
\[d_{k}(\la) = \sum_{(i,j) \in \la} {(c_{ij}^{\alpha})}^k.\]
For any integers $n,p,k \ge 0$ we set
\[F_{npk}(\la) = P_{npk}(d_{1}(\la),d_{2}(\la),d_{3}(\la),\ldots)
=\sum_{|\mu| = n} \frac{\npbi{\mu}{p}{k}}{z_{\mu}}\, d_\mu(\la)\]
and $F_{nk}(\la)=F_{n0k}(\la)$, with the obvious notation 
\[d_\mu(\la)= \prod_{i \ge 1} d_i(\la)^{m_i(\mu)}.\]
In other words we choose the following specialization 
\begin{equation*}
 X_0=d_{0}(\la)=|\la| \quad,\quad X_{k}=d_{k}(\la) \quad  (k \ge 1),
\end{equation*}
though the quantities $d_{k}(\la)$ are no longer independent indeterminates.

In~\cite{LL} we have shown the following development
\begin{equation}
\frac{(x+y)_{\la}}{(y)_{\la}} =
\sum_{i,j\ge 0} (-1)^{j} \frac{x^i}{y^{i+j}}
\left(\sum_{k=0}^{\mathrm{min}(i,j)}\binom{|\la|-j}{i-k}\,F_{jk}(\la) \right),
\end{equation}
where the sum takes over any  $j\ge 0$ and not only over $|\la|-j\ge 0$. 

An analogous result can be deduced from Theorem 3.1.
\begin{theorem}
Let $x,y$ be two indeterminates. Then we have
\begin{multline*}
\frac{(x+y+1)_{\la}}{(x+y)_{\la}}\, \frac{(x)_{\la}}{(x+1)_{\la}} =
\sum_{n\ge 0} {\Big(-\frac{y}{x^2}\Big)}^n 
{\Big(1+\frac{y+1}{x}\Big)}^{-n}\\
\sum_{p,q \ge 0} {\Big(\frac{-1}{x}\Big)}^{p+q} 
{\Big(1+\frac{y+1}{x}\Big)}^{-q}\left( \sum_{k=0}^{\mathrm{min}(n,p+q)}
\binom{|\la|+n-1}{n-k} F_{p+q,p,k}(\la) \right).
\end{multline*}
\end{theorem}
\begin{proof}
With $c_{ij}=c_{ij}^{\alpha}$ for short, the left-hand side can be written
\begin{equation*}
\begin{split}
\textrm{LHS}&=\prod_{(i,j)\in \la} \frac{(x+y+1+c_{ij})}{(x+y+c_{ij})}
\, \frac{(x+c_{ij})}{(x+1+c_{ij})}\\
&= \prod_{(i,j)\in \la} {\Big(1+\frac{y}{(x+y+1+c_{ij})
(x+c_{ij})}\Big)}^{-1}.
\end{split}
\end{equation*}
Setting $u=-1/x$, $v=-1/(x+y+1)$, and using the classical series expansion
\[\textrm{log}(1-a)=-\sum_{k\ge1}\frac{a^k}{k},\]
we get
\begin{equation*}
\begin{split}
\textrm{LHS}&=\prod_{(i,j)\in \la} \textrm{exp} \Big[
\sum_{k\ge1}\frac{(-y)^k}{k}{\Big((x+y+1+c_{ij})(x+c_{ij})\Big)}^{-k}\Big]\\
&= \textrm{exp}  \Big[\sum_{(i,j)\in \la} \Big(\sum_{k\ge1} 
\frac{(-yuv)^k}{k} {(1-c_{ij} u)}^{-k} {(1-c_{ij} v)}^{-k}\Big)\Big].
\end{split}
\end{equation*}
Then using the series expansion
\[{(1-a)}^{-k}=\sum_{r\ge0} a^r \frac{(k)_r}{r!},\]
we obtain
\begin{equation*}
\begin{split}
\textrm{LHS}&= \textrm{exp}  \Big[ \sum_{k\ge1} 
\frac{(-yuv)^k}{k} \Big( \sum_{r,s\ge0}  u^r v^s \frac{(k)_r}{r!}  \frac{(k)_s}{s!}
d_{r+s}(\la)\Big)\Big]\\
&= \prod_{k\ge1} \Big( \sum_{m_k\ge0} \frac{1}{m_k!} {\Big[ \frac{(-yuv)^k}{k} 
\Big( \sum_{r,s\ge0}  u^r v^s \frac{(k)_r}{r!} \,
\frac{(k)_s}{s!} \, d_{r+s}(\la) \Big) \Big]}^{m_k} \Big)\\
&=\sum_\mu \frac{(-yuv)^{|\mu|}}{z_{\mu}}
\prod_{k\ge1} {\Big(\sum_{r,s\ge0}  u^r v^s \frac{(k)_r}{r!}\,  
\frac{(k)_s}{s!}\, d_{r+s}(\la)\Big)}^{m_k(\mu)}.
\end{split}
\end{equation*}
We conclude by applying Theorem 3.1 with the specialization 
$X_{k}=d_{k}(\la)$, which gives
\[\textrm{LHS}=\sum_{n\ge 0} (-yuv)^n
\sum_{p,q \ge 0}  u^p v^q
\left( \sum_{k=0}^{\mathrm{min}(n,p+q)} \binom{|\la|+n-1}{n-k} F_{p+q,p,k}(\la) 
\right).\]
\end{proof}

\begin{coro}
Let $x,y$ be two indeterminates. Then we have
\[\frac{(x+y+1)_{\la}}{(x+y)_{\la}}\, \frac{(x)_{\la}}{(x+1)_{\la}} =
\sum_{r\ge 0} c_r{\left(\frac{-1}{x}\right)}^r\]
with
\[c_{r}=\sum_{\begin{subarray}{c}n,p,q \ge 0 \\ 2n+p+q \le r\end{subarray}} 
(-y)^n (y+1)^p\binom{n+p+q-1}{p}
\left(\sum_{k=0}^{\mathrm{min}(n,r-2n-p)}
\binom{|\la|+n-1}{n-k} F_{r-2n-p,q,k}(\la) \right).\]
\end{coro}

\begin{remark} 
Observe the situation for low indices. One 
has $c_0 =1$ since $F_{000}(\la)=1$, and $c_1=0$ since $F_{1q0}(\la)=0$.
\end{remark}

This result is useful in the study of Jack 
polynomials, allowing to prove some conjectures stated 
in~\cite{La4,La5}. The proofs will be given in Sections 8 and 9.

\section{Symmetric functions}

Let $A=\{a_1,a_2,a_3,\ldots\}$ a (possibly infinite) set of
independent indeterminates ($A$ is called an alphabet). The generating functions
\begin{equation*}
\begin{split}
E_t(A) &= \prod_{a\in A} (1 +ta) =\sum_{k\geq0} t^k\, e_k(A)\\
H_t(A) &=  \prod_{a\in A}  \frac{1}{1-ta} = \sum_{k\geq0} t^k\, h_k(A)\\
P_t(A) &= \sum_{a\in A} \frac{a}{1-ta} =\sum_{k\geq1} t^{k-1}p_k(A) 
\end{split}
\end{equation*}
define symmetric functions known as respectively elementary, complete and 
power sums. For any partition $\mu$, we define functions $e_\mu$, 
$h_\mu$ or $p_\mu$ by
\[f_{\mu}= \prod_{i=1}^{l(\mu)}f_{\mu_{i}}=\prod_{k\geq1}f_k^{m_{k}(\mu)},\]
where $f_{i}$ stands for $e_i$, $h_i$  or $p_i$. The monomial symmetric function 
$m_\mu(A)$ is defined as the sum of all {\it distincts} monomials 
$\prod_{i}a_i^{m_i}$ such that $(m_i)$ is a permutation of $\mu$~\cite{Ma}.

When $A$ is infinite, each of the three sets of functions $e_i(A)$, 
$h_i(A)$ or $p_i(A)$ forms an algebraic basis of
$\mathbf{S}[A]$, the symmetric functions algebra of $A$.
Each of the sets of functions $e_\mu(A)$, $h_\mu(A)$, $p_\mu(A)$ is a linear basis of this 
algebra. It is thus possible to define the symmetric 
algebra $\mathbf{S}$ abstractly, as the 
$\mathbf{R}$-algebra generated by the functions $e_i$, $h_i$ or $p_i$.

This is no longer true when $A$ is finite, a situation which is often 
encountered. In that case the functions $e_\mu(A)$ (resp. $h_\mu(A)$, 
$p_\mu(A)$) are no longer linearly independent. We have encountered 
such a situation with the previous specialization $X_k=d_k(\la)$. 
Let 
\[A_\la=\left\{c_{ij}^{\alpha}=j-1-(i-1)/\alpha,\quad (i,j) \in \la \right\}\]
denote the finite alphabet of $\alpha$-contents of $\la$.
By definition we have
\[d_k(\la)=p_k(A_\la) \quad  (k \ge 1),\quad d_\mu(\la)=p_\mu(A_\la).\]

Thus the left-hand sides of the 
relation of Theorem 5.1, as well as of relation (5.1), 
are symmetric functions of the alphabet $A_\la$. The 
expressions given are in terms 
of the power-sums $p_\mu(A_\la)$, but it might be also interesting to 
convert them in terms of any other symmetric functions of $A_\la$.

With this conversion in mind, we are led to the two following open problems.
Let $A=\{a_1,a_2,a_3,\ldots\}$ be any (finite or infinite) alphabet. 
Choose the specialization  $X_k = p_k(A), \, (k\ge 1)$. By 
definition for any $n,k \ge 0$ we have
\[P_{nk}(-X) = 
\sum_{|\mu| = n} (-1)^{l(\mu)}
\frac{\pbi{\mu}{k}}{z_\mu}\, p_\mu(A).\]
In~\cite{LL} (Lemma 2, p. 306) we have shown 
\[P_{nk}(-X) =
(-1)^{k}\sum_{\begin{subarray}{c}|\mu| = n\\
l(\mu) = k\end{subarray}} m_{\mu} (A).\]
Let $\omega$ be the involution defined by $\omega(p_r)= (-1)^{r-1} p_r$. 
One has $\omega(P_{nk}(X)) =(-1)^n P_{nk}(-X)$, hence
\[P_{nk}(X) = (-1)^{n-k}\sum_{\begin{subarray}{c}|\mu| = n\\
l(\mu) = k\end{subarray}} f_{\mu} (A),\]
with $f_{\mu}$ the fourth classical basis of ``forgotten''symmetric 
functions (\cite{Ma}, p. 22).

\begin{prob}
Evaluate $P_{npk}(-X)$ in terms of monomial symmetric 
functions of $A$.
\end{prob}

Using ACE~\cite{V} there is some experimental evidence that
\[P_{npk}(-X) =
(-1)^{k} \sum_{\begin{subarray}{c}|\mu| = n\\
l(\mu) = k\end{subarray}} \chi_{\mu}\,m_{\mu} (A),\]
with $\chi_{\mu}$ some polynomial in the multiplicities 
of $\mu$. It seems that for $p \le 3$ one has
\[\chi_{\mu}= \binom{k+p-1}{p}-\binom{k+p-3}{p-2} m_1(\mu)-
\binom{k+p-4}{p-3} m_2(\mu).\]

Now if we define the alphabet 
\[B_0=\{a_1/(1-a_1),a_2/(1-a_2),a_3/(1-a_3),\ldots\},\] 
it was proved in~\cite{JZ} (Theorem 2) that for 
any $k\ge 1$, one has
\begin{equation*}
\begin{split}
p_k(B_0) = \sum_{n\ge k} &\binom{n-1}{k-1} p_{n}(A)\\
h_k(B_0) = \sum_{n\ge k}P_{nk}(X) \quad&,\quad
e_k(B_0) = (-1)^k \sum_{n\ge k}P_{nk}(-X). 
\end{split}
\end{equation*}

\begin{prob}
For any integer $p\ge 0$, find an alphabet $B_p$ such that for 
any $k\ge 1$,
\[p_k(B_p) = \sum_{n\ge \mathrm{max}(k,p)} \frac{k}{n} \, \nbi{n}{p}{k} p_{n}(A).\]
Give the expansion of $h_k(B_p)$ (resp. $e_k(B_p)$) in terms of the quantities $P_{npk}(X)$ 
(resp. $P_{npk}(-X)$).
\end{prob}

For any two alphabets $A$ and $B$, their difference $A-B$ 
(which is not their difference as sets) is defined by
\[E_t(A-B) = E_t(A)\,{E_t(B)}^{-1} \quad,\quad H_t(A-B) = 
H_t(A)\,{H_t(B)}^{-1}.\]
In Sections 8 and 9 we shall need the following result about Lagrange interpolation~\cite{Ls}.
\begin{laglem}
Let $A$ and $B$ be two finite alphabets with respective 
cardinals $n$ and $m$. For any integer $r\ge0$ we have
\[\sum_{a \in A} a^r \, \frac{\displaystyle\prod_{b \in B}(a-b)}
{\displaystyle\prod_{c \in A, \, c \neq a} (a-c)} =
h_{m-n+r+1}(A-B).\]
\end{laglem}
Alain Lascoux~\cite{Ls} mentions that when $B$ is empty, 
this result was already known to Euler.

\section{Shifted symmetric functions}

Though the theory of symmetric functions goes back to the early 19th century, the 
notion of ``shifted symmetric'' functions is very recent. We refer 
to~\cite{K,KS,Sa}, \cite{Ok1,Ok2,Ok3} and to references given there. We shall follow 
the presentation given by \cite{Ok1,Ok2,Ok3}. 

Let $\alpha$ be some fixed positive real number and $\mathbf{F}$ be the field of rational 
functions in $\alpha$. A polynomial 
in $n$ variables $(x_1,x_2,\ldots,x_n)$ with coefficients in 
$\mathbf{F}$ is said to be shifted symmetric if it is symmetric 
in the ``shifted variables'' 
$x_i-i/\alpha$. Let ${\mathbf{S}}_n^{*}$ denote this subalgebra of 
$\mathbf{F}[x_1,x_2,\ldots,x_n]$.

Now let $X=\{x_1,x_2,x_3,\ldots\}$ be an infinite alphabet. Consider 
the morphism from ${\mathbf{S}}_{n+1}^{*}$ to 
${\mathbf{S}}_n^{*}$ given by $x_{n+1}=0$. In analogy with symmetric functions, the shifted 
symmetric algebra $\mathbf{S}^{*}[X]$ is defined as the projective limit of the
algebras ${\mathbf{S}}_n^{*}$ with respect to this morphism. In other words 
a ``shifted symmetric function'' $f\in \mathbf{S}^{*}[X]$ is 
a family $\{f_n, n\ge 1\}$ with the two following properties
\begin{enumerate}
    \item[(i)]  $f_n \in {\mathbf{S}}_n^{*}$ (shifted symmetry),
    \item[(ii)]  $f_{n+1}(x_1,x_2,\ldots,x_n,0)=f_{n}(x_1,x_2,\ldots,x_n)$ 
    (stability).
\end{enumerate}
For instance if for any integer $k\ge 1$ we define the ``shifted 
power sums'' by
\[p_k^{*}(x)=\sum_{i\ge 
1}\Big([x_i-(i-1)/\alpha]_k-[-(i-1)/\alpha]_k\Big),\]
these polynomials algebraically generate ${\mathbf{S}}^{*}[X]$.

Any element $f\in \mathbf{S}^{*}[X]$ may be evaluated at any sequence 
$x=(x_1,x_2,\ldots)$ with finitely many non zero terms, hence at any partition 
$\la$. Moreover by analyticity, $f$ is entirely determined by its restriction 
$f(\la)$ to partitions. From now on we shall perform this
identification, and consider $\mathbf{S}^{*}[X]$ as a function 
algebra $\mathbf{S}^{*}[\mathcal{P}]$ on the set $\mathcal{P}$ of partitions. 
In this approach, a partition is no longer a parameter but an argument.
The following lemma is crucial (\cite{La4}, p. 151). The argument is taken 
from~\cite{Ok4} (Lemma 7.1).

\begin{prop}
For any integer $k\ge 1$, the function $\la \rightarrow d_k(\la)$ defines a shifted 
symmetric function.
\end{prop}
\begin{proof}
Let $t(k,m)$ denote the inverse matrix of $s(k,m)$, the matrix of Stirling numbers of 
the first kind, that is $x^k=\sum_{m=1}^kt(k,m)[x]_{m}$. 
We have
\[d_k= \sum_{m=1}^k t(k,m) \frac{p_{m+1}^{*}}{m+1}.\]
Indeed we can write
\begin{equation*}
\begin{split}
d_k(\la)&= \sum_{m=1}^k \sum_{(i,j)\in \la} t(k,m)\,
[j-1-(i-1)/\alpha]_{m}\\
&= \sum_{m=1}^k \frac{t(k,m)}{m+1} \sum_{i= 1}^{l(\la)} 
\Big([\la_i-(i-1)/\alpha]_{m+1}-[-(i-1)/\alpha]_{m+1}\Big),
\end{split}
\end{equation*}
the last equation being a direct consequence of the identity
$m[x]_{m-1}=[x+1]_{m}-[x]_{m}$.
\end{proof}

For instance we have
\[d_1=\frac{p_2^{*}}{2} \quad,\quad d_2=
\frac{p_3^{*}}{3}+\frac{p_2^{*}}{2} \quad,\quad
d_3= \frac{p_4^{*}}{4}+p_3^{*}+\frac{p_2^{*}}{2}.\]

For each element $f\in \mathbf{S}^{*}[X]$ we define 
its so called leading symmetric term $[f]\in \mathbf{S}[X]$ 
as the highest degree homogeneous 
part of $f$. For instance we have $[p_k^{\star}]=p_k$. The map $f\rightarrow [f]$ 
provides an isomorphism of the graded 
algebra associated to the filtered algebra $\mathbf{S}^{*}[X]$ onto the symmetric function 
algebra $\mathbf{S}[X]$.
Assuming that the leading terms $[f_1],[f_2],\ldots,[f_n]$ of a sequence 
$f_1,f_2,\ldots,f_n$ generate 
$\mathbf{S}[X]$, then this sequence itself generates 
$\mathbf{S}^{*}[X]$. 

Now for $k\ge1$ as a consequence of Lemma 7.1,
\[[d_k]=\frac{[p_{k+1}^{*}]}{k+1} =\frac{p_{k+1}}{k+1}.\]
Thus the shifted symmetric functions $p_1^{*}=p_1$ and 
$d_k, k\ge 1$ algebraically generate $\mathbf{S}^{*}[X]$.
In other words any $f\in \mathbf{S}^{*}[\mathcal{P}]$ 
may be expanded in terms of the functions $|\la|^r d_\mu(\la)$ with $r\ge0$ 
and $\mu$ any partition. 

But as seen before $d_\mu(\la)=p_\mu(A_\la)$ and there is no 
reason that we shall restrict to the basis of power sums in the alphabet $A_\la$. 
Any other basis of the symmetric algebra may be convenient. Therefore
any $f\in \mathbf{S}^{*}[\mathcal{P}]$ may be expanded in terms of
the functions $|\la|^r b_\mu(A_\la)$, with $r\ge0$ and $b_\mu$  
any linear basis of the symmetric algebra $\mathbf{S}$, viewed 
abstractly.

Finally we have proved
\begin{theorem}
Let $G=\mathbf{F}[w]$ be the field of polynomials in some 
indeterminate $w$ with coefficients 
in $\mathbf{F}$. The evaluation map $p\otimes f \rightarrow p(|\la|)f(A_\la)$ 
is an isomorphism of $G\otimes\mathbf{S}$ onto $\mathbf{S}^{*}[\mathcal{P}]$.
\end{theorem}

Let us give a simple but enlightening example. Recall that 
the generalized raising 
and lowering factorial are given by
\[(x)_{\la}=\prod_{a \in A_{\la}} (x+a) \quad,\quad 
[x]_{\la}=\prod_{a \in A_{\la}} (x-a).\] 
Then let us define (\cite{La4}, p. 64)
\[(x)_{\la}= \sum_{k\geq 0} c_k(\la) \, x^{|\la|-k} \quad , \quad
\frac{1}{[x]_{\la}}= \sum_{k\geq 0} C_k(\la) \, \frac{1}{x^{|\la|+k}}.\]
The quantities $c_k(\la)$ and $C_k(\la)$ are 
generalizations of Stirling numbers 
of the first and second kind, since when $\la$ is a row-partition $(n)$, 
they are respectively $|s(n,k)|$ and $S(n+k-1,n-1)$.

We know immediatly that $c_k(\la)$ and $C_k(\la)$ are 
shifted symmetric functions, i.e. belong to 
$\mathbf{S}^{*}[\mathcal{P}]$, and that they satisfy
\[c_k(\la)=\sum_{|\mu|=k} (-1)^{k-l(\mu)} \frac{1}{z_\mu} \, d_\mu(\la)
\quad , \quad
C_k(\la)=\sum_{|\mu|=k} \frac{1}{z_\mu} \, d_\mu(\la).\]
Indeed their definition are merely the generating functions 
for elementary and complete symmetric functions of the 
alphabet $A_{\la}$ since
\begin{equation*}
\begin{split}
(x)_{\la}=\prod_{a \in A_{\la}} (x+a)
&= x^{|\la|} E_{1/x}(A_\la) = \sum_{k\geq0} x^{|\la|-k}\, e_k(A_\la)\\
\frac{1}{[x]_{\la}}=\prod_{a \in A_{\la}} \frac{1}{x-a} &=\frac{1}{x^{|\la|}} 
H_{1/x}(A_\la) = \sum_{k\geq0} \frac{1}{x^{|\la|+k}}\, h_k(A_\la).
\end{split}
\end{equation*}
Thus it is enough to
apply the classical Cauchy formulas (\cite{Ma}, p. 25) to get
\begin{equation*}
\begin{split}
c_k(\la)=e_k(A_\la)&=\sum_{|\mu|=k} (-1)^{k-l(\mu)}
\frac{1}{z_\mu} \, p_\mu(A_\la)\\
C_k(\la)=h_k(A_\la)&=\sum_{|\mu|=k} \frac{1}{z_\mu} \, p_\mu(A_\la).
\end{split}
\end{equation*}

\section{Application to Jack polynomials}

The reference for Jack polynomials is Chapter 6, Section 10 of the book of 
Macdonald \cite{Ma}.
Let $\alpha$ be some fixed positive real number. Jack 
polynomials $J_\la^{\alpha}$ form a basis of $\mathbf{F}\otimes\mathbf{S}$, the algebra of 
symmetric functions with coefficients in $\mathbf{F}$. 

Jack polynomials 
satisfy the following generalization of Pieri formula~\cite{St}. 
For any partition $\la$ and any integer $i$ such that $1 \le i \le l(\la)+1$, 
denote $\la^{(i)}$ the partition $\mu$, if it exists, such that 
$\mu_j=\la_j$ for $j\neq i$ and $\mu_i=\la_i +1$. Then we have
\[e_1 \,J_\la^{\alpha} =\sum_{i=1}^{l(\la)+1} 
c_i^{\alpha}(\la) \, J_{\la^{(i)}}^{\alpha}.\]
The Pieri coefficients $c_i^{\alpha}(\la)$ have the following analytic 
expression~\cite{La8}
\[c_i^{\alpha}(\la) = \frac {1}{\alpha \la_i+l(\la)-i+2}
\prod_{\begin{subarray}{c}j=1 \\ j \neq i\end{subarray}}^{l(\la)+1} 
\frac{\alpha(\la_i-\la_j)+j-i+1}
{\alpha(\la_i-\la_j)+j-i}.\]

For any integer $r\ge 0$ we define
\[s_r(\la) =\sum_{i=1}^{l(\la)+1} 
{\left(\la_i-\frac{i-1}{\alpha}\right)}^r c_i^{\alpha}(\la).\]
In~\cite{La4} we gave the values of $s_r(\la)$ up to 
$r=9$. For instance one has
\[s_0(\la) = 1 \quad,\quad s_1(\la) = 0 \quad,\quad
s_2(\la)=|\la|/\alpha \quad,\quad s_3(\la)=2d_1(\la)/\alpha + 
|\la|(\alpha-1)/\alpha^2.\] 

Since $c_i^{\alpha}(\la)$ is a rational function of the variables 
$\{\la_i-i/\alpha\}$, it is not obvious that $s_r(\la)$ is a shifted 
symmetric function, i.e. a \textit{polynomial} in $\{\la_i-i/\alpha\}$. 
A proof has been given by Macdonald (\cite{La4}, Theorem 6.1, p. 69). 
A more direct proof is given below, together with a new explicit formula. 

In~\cite{La4} (Conjecture 6.2, p. 70) we had stated the following 
conjecture.
\begin{conj}
For any partition $\la$ and any integer $r\ge 0$ one has
\[s_r(\la) =\sum_{i= 0}^{[r/2]} \, \sum_{j=0}^{r-2i} \, 
\sum_{k=0}^{\mathrm{min}(i,j)} \frac{1}{\alpha^i}
\left(1-\frac{1}{\alpha}\right)^{r-2i-j}
\binom{|\la|+i-1}{i-k}
\sum_{|\rho|=j} u_{ijk}^{\rho} (r) \frac{d_\rho(\la)}{z_{\rho}},\]
where the coefficients $u_{ijk}^{\rho} (r)$ are positive integers.
\end{conj} 
\begin{theorem}
The previous conjecture is true, with
\[u_{ijk}^{\rho} (r) = \sum_{s=0}^j \npbi{\rho}{s}{k} \binom{r+s-i-j-1}{r-2i-j}.\]
In other words we have
\begin{multline*}
s_r(\la)=\sum_{\begin{subarray}{c}n,p,q \ge 0 \\ 2n+p+q \le r\end{subarray}}  
\frac{1}{\alpha^n} 
\left(1-\frac{1}{\alpha}\right)^p \binom{n+p+q-1}{p}\\
\left(\sum_{k=0}^{\mathrm{min}(n,r-2n-p)}
\binom{|\la|+n-1}{n-k} F_{r-2n-p,q,k}(\la) \right).
\end{multline*}
\end{theorem}
\begin{proof}
We shall prove that the quantities $s_r(\la)$ have the following generating 
function
\[\sum_{r\ge 0} s_r(\la){\left(\frac{-1}{x}\right)}^r
=\frac{(x-1/\alpha+1)_{\la}}{(x-1/\alpha)_{\la}}
\frac{(x)_{\la}}{(x+1)_{\la}}.\]
Then the statement will immediatly follow from Corollary 5.2 with
$ y=-1/\alpha$.
To prove this generating function we apply Lagrange Lemma for the two following alphabets
\begin{equation*}
\begin{split}
A=\{a_i =\alpha\la_i -i+1 & ,\quad i=1,\ldots,l(\la)+1 \},\\
B=\{b_i =\alpha\la_i -i & ,\quad i=1,\ldots,l(\la)\}.
\end{split}
\end{equation*}
Then it is obvious that
\[c_i^{\alpha}(\la)= \frac {\displaystyle \prod_{b \in B}(a_i-b)}
{\displaystyle \prod_{c \in A,\, c \neq a_i}(a_i-c)},\]
and the Lagrange Lemma yields
\[\alpha^r s_r(\la)= h_r(A-B).\]
The generating function of $\alpha^r s_r(\la)$ is thus
\begin{equation*}
\begin{split}
H_z(A-B)&= \frac {\prod_{b\in B} (1-zb)}{\prod_{a\in A} (1-za)}\\
&=\prod_{i=1}^{l(\la)} \frac{1-z(\alpha\la_i -i)}{1-z(\alpha\la_i-i+1)} \, 
\frac{1}{1+zl(\la)}.
\end{split}
\end{equation*}
With $\alpha z=-1/x$ we obtain
\[\sum_{r\ge 0} s_r(\la){\left(\frac{-1}{x}\right)}^r =
\prod_{i=1}^{l(\la)} \frac{x+\la_i -i/\alpha}{x+\la_i -(i-1)/\alpha} \, 
\frac{x}{x-l(\la)/\alpha}.\]
Since
\[\frac {(u+1)_\la}{(u)_\la}=\prod_{i=1}^{l(\la)}
\prod_{j=1}^{\la_i} \frac{u+j-(i-1)/\alpha}{u+j-1-(i-1)/\alpha}
=\prod_{i=1}^{l(\la)} \frac{u+\la_i-(i-1)/\alpha}{u-(i-1)/\alpha}\]
we conclude easily.
\end{proof}

\begin{remark} 
As already emphasized, Theorem 8.1 gives the 
expansion of the shifted symmetric function $s_r(\la)$ 
over the basis $|\la|^r p_\mu(A_\la)$. 
An open problem is to find explicit expansions over different 
basis of the 
symmetric algebra $\mathbf{S}[A_{\la}]$, for instance  $|\la|^r m_\mu(A_\la)$.
\end{remark}

\section{Application to shifted Jack polynomials}

For any partition $\mu$ there exists a shifted symmetric function 
$\widetilde{J_\mu^{\alpha}}$ such that
\begin{enumerate}
\item[(i)]  degree $\widetilde{J_\mu^{\alpha}} = |\mu|$,
\item[(ii)] $\widetilde{J_\mu^{\alpha}}(\la)=0$ except if 
$\mu_i \le \la_i$ for any $i$, and 
$\widetilde{J_\mu^{\alpha}}(\mu) \neq 0$.
\end{enumerate}

This function is called the shifted Jack polynomial associated with 
$\mu$ ~\cite{K,KS,Ok1,Ok2,Ok3,Sa}. 
It is unique up to the value of $\widetilde{J_\mu^{\alpha}}(\mu)$.  
When conveniently normalized, $\widetilde{J_\mu^{\alpha}}(\la)$ is merely the 
generalized binomial coefficient
\[\frac{\widetilde{J_\mu^{\alpha}}(\la)}{\widetilde{J_\mu^{\alpha}}(\mu)} 
= {\binom{\la}{\mu}}_{\alpha}\]
appearing in the generalized binomial formula for Jack polynomials
\cite{Ka,La6,La7,Ok2}. 

We mention the two following remarkable facts :
\begin{enumerate}
\item[(a)] In the definition of $\widetilde{J_\mu^{\alpha}}$, 
the overdetermined system of linear conditions (ii) may be replaced by
the weaker conditions

\item[(iii)] $\widetilde{J_\mu^{\alpha}}(\la)=0$ except if $|\mu| \le |\la|$, 
and $\widetilde{J_\mu^{\alpha}}(\mu) \neq 0$.
\item[(b)] For some constant $h_\mu$ one has
$[\widetilde{J_\mu^{\alpha}}]= h_\mu J_\mu^{\alpha}$.
\end{enumerate}

For any partition $\la$ and any integer $i$ such that $1 \le i \le l(\la)$, 
we denote $\la_{(i)}$ the partition $\mu$, if it exists, such that 
$\mu_j=\la_j$ for $j\neq i$ and $\mu_i=\la_i -1$. 
In \cite{La6}, p. 320 (see also \cite{La7}, Theorem 5), we have proved 
that
\[{\binom{\la}{\la_{(i)}}}_{\alpha} = \left(\la_i +\frac {l(\la)-i}{\alpha}\right)
\prod_{\begin{subarray}{c}j=1 \\ j \neq i \end{subarray}}^{l(\la)} 
\frac{\alpha(\la_i-\la_j)+j-i-1}
{\alpha(\la_i-\la_j)+j-i}.\]
For any integer $r\ge 0$ we introduce the shifted symmetric function
\[\sigma_r(\la) =\sum_{i=1}^{l(\la)} {\left(\la_i-\frac{i-1}{\alpha}\right)}^r 
{\binom{\la}{\la_{(i)}}}_{\alpha}.\]
As shown in \cite{La6}, Theorem 8 (see also \cite{La7}, Theorem 3) 
we have $\sigma_0(\la)= |\la|$. But even the computation of the first values
\[\sigma_1(\la) = 2d_1(\la) + |\la|\quad,\quad \sigma_2(\la) = 3d_2(\la)
+\left(3+\frac{1}{\alpha}\right)d_1(\la) + |\la|-\frac{1}{\alpha} \binom{|\la|}{2} \]
seems rather difficult, except by the method given below.

\begin{theorem}
We have
\[\sigma_r(\la)=c_{r+1}-\alpha \, c_{r+2},\]
with
\begin{multline*}
c_{r}=\sum_{\begin{subarray}{c}n,p,q \ge 0 \\ 2n+p+q \le r\end{subarray}}
\left(-\frac{1}{\alpha}\right)^n
\left(1+\frac{1}{\alpha}\right)^p \binom{n+p+q-1}{p}\\
\left(\sum_{k=0}^{\mathrm{min}(n,r-2n-p)}
\binom{|\la|+n-1}{n-k} F_{r-2n-p,q,k}(\la) \right).
\end{multline*}
\end{theorem}
\begin{proof}
We shall prove that the quantities $\sigma_r(\la)$ have the following generating 
function
\[\sum_{r\ge 0} \sigma_r(\la){\left(\frac{-1}{x}\right)}^r =
-x(\alpha x+1)\, \left(\frac{(x+1/\alpha+1)_{\la}}{(x+1/\alpha)_{\la}}
\frac{(x)_{\la}}{(x+1)_{\la}} -1\right)\]
Then the statement will immediatly follow from Corollary 5.2 with
$ y=1/\alpha$.
To prove this generating function we apply Lagrange Lemma for the two following alphabets
\begin{equation*}
\begin{split}
A=\{a_i =\alpha\la_i -i+1 & ,\quad i=1,\ldots,l(\la) \},\\
B=\{b_i =\alpha\la_i -i+2 & ,\quad i=1,\ldots,l(\la)+1\}.
\end{split}
\end{equation*}
Then it is obvious that
\[{\binom{\la}{\la_{(i)}}}_{\alpha}= 
-\frac{1}{\alpha}\, \frac {\displaystyle \prod_{b \in B}(a_i-b)}
{\displaystyle \prod_{c \in A,\, c \neq a_i}(a_i-c)},\]
and the Lagrange Lemma yields
\[-\alpha^{r+1} \sigma_r(\la)= h_{r+2}(A-B).\]
The generating function of $-\alpha^{r+1} \sigma_r(\la)$ is thus
\[-\sum_{r\ge 0} \alpha^{r+1} z^r \sigma_r(\la) =
\frac{1}{z^2} \sum_{r\ge2}z^r h_{r}(A-B) = 
\frac{1}{z^2}\big(H_z(A-B)-zh_{1}(A-B)-1\big).\]
But we have
\begin{equation*}
\begin{split}
H_z(A-B)&= \frac {\prod_{b\in B} (1-zb)}{\prod_{a\in A} (1-za)}\\
&=\prod_{i=1}^{l(\la)} \frac{1-z(\alpha\la_i -i+2)}{1-z(\alpha\la_i-i+1)} \, 
(1+z(l(\la)-1).
\end{split}
\end{equation*}
And it is easily checked that
\[h_{1}(A-B)=
\sum_{i=1}^{l(\la)} (a_i -b_i)+l(\la)-1=-1.\]
With $\alpha z=-1/x$ we obtain
\begin{equation*}
\sum_{r\ge 0} \sigma_r(\la){\left(\frac{-1}{x}\right)}^r =
-\alpha x^2 \left(
\prod_{i=1}^{l(\la)} \frac{x+\la_i -(i-2)/\alpha}{x+\la_i -(i-1)/\alpha} \, 
\frac{x-(l(\la)-1)/\alpha}{x} -\frac{1}{\alpha x} -1 \right).
\end{equation*}
The right-hand side can be written as
\[-\alpha x \, (x-(l(\la)-1)/\alpha) \,
\frac{(x+1/\alpha+1)_{\la}}{(x+1/\alpha)_{\la}}
\frac{(x)_{\la}}{(x+1)_{\la}}\,\prod_{i=1}^{l(\la)} 
\frac{x -(i-2)/\alpha}{x -(i-1)/\alpha} 
+x +\alpha x^2.\]
Hence the result.
\end{proof}

\section{Probability distributions}

Let us briefly outline the connections of these 
results with probability distributions on the set 
of Young diagrams~\cite{Ke,Ke2,OR}.

We denote by $\mathcal{P}_n$ the set of Young (i.e. Ferrers) 
diagrams with $n$ cells, and by $\mathcal{P}=\cup\,  
\mathcal{P}_n$ the lattice of all Young diagrams ordered by 
inclusion. We write $\la \nearrow \Lambda$ if the diagram $\Lambda$ 
is obtained from $\la$ by adding a cell, i.e. 
$\Lambda=\la^{(i)}$ for some row $i$. We consider such a 
pair $(\la,\Lambda)$ as an ``oriented edge'' of the ``infinite 
graph'' $\mathcal{P}$.

We use the Pieri formula 
\[e_1 \,J_\la^{\alpha} =\sum_{\Lambda :\la \nearrow \Lambda} 
\kappa(\la,\Lambda) \, J_{\Lambda}^{\alpha}.\]
to define a ``multiplicity function'' $\kappa$ 
on the set of edges of $\mathcal{P}$. 
A function $\Phi$ defined on $\mathcal{P}$ is said to be harmonic if
for all vertices $\la \in\mathcal{P}$ one has
\[\Phi(\la)=\sum_{\Lambda :\la \nearrow \Lambda}
\kappa(\la,\Lambda) \, \Phi(\Lambda).\]
Let $\varphi$ be the particular case corresponding to 
the initial condition $\varphi(\emptyset) =1$. Then (\cite{Ke}, 
Lemma 7.2) one has $\varphi(\la) =1$ 
for any $\la \in\mathcal{P}$.

Similarly we define a ``dimension function'' $\mathrm{dim}\ \Lambda$ 
by the following recurrence
\[\mathrm{dim}\ \Lambda=\sum_{\la :\la \nearrow \Lambda}
\kappa(\la,\Lambda) \, \mathrm{dim}\ \la\]
with the initial condition 
$\mathrm{dim}\ \emptyset =1$. Then (\cite{Ke}, 
Corollary 6.10) one has 
\[\mathrm{dim}\ \la =|\la| ! \, \frac{\alpha^{|\la|}}{j_\la^{\alpha}}\]
with $j_\la^{\alpha}=\langle J_\la^{\alpha},J_\la^{\alpha} \rangle$ 
the norm of $J_\la^{\alpha}$ (\cite{St}, p. 97).

It is shown in \cite{Ke,Ke2} that for any oriented edge $\la \nearrow \Lambda$ 
the quantities
\[p_{\la}(\Lambda)=\kappa(\la,\Lambda)\,
\frac{\varphi(\Lambda)}{\varphi(\la)}
\quad,\quad
q_{\Lambda}(\la)=\kappa(\la,\Lambda)\, 
\frac{\mathrm{dim}\ \la}{\mathrm{dim}\ \Lambda}\]
define two discrete ``transition'' 
$p_{\la}$ and ``co-transition'' $q_{\Lambda}$ probability 
distributions, associated with the Young diagrams $\la$ and 
$\Lambda$ respectively.

But as a consequence of \cite{La7}, Theorem 2 
for any oriented edge $\la \nearrow \Lambda$ one has
\[\kappa(\la,\Lambda)=\alpha \ \frac{j_\la^{\alpha}}{j_\Lambda^{\alpha}}
\ {\binom{\Lambda}{\la}}_{\alpha},\]
so that with our previous notations, we can write
\[p_{\la}(\la^{(i)})=c_i^{\alpha}(\la)\quad,\quad 
q_{\Lambda}(\Lambda_{(i)})=\frac{1}{|\Lambda|}
{\binom{\Lambda}{\Lambda_{(i)}}}_{\alpha}.\]
Now let us consider the $\alpha$-content random variable 
$c$ defined on any cell $(i,j)$ by 
$c(i,j)=j-1-(i-1)/\alpha$. If $\la \nearrow \Lambda$ 
is an oriented edge 
(i.e. $\Lambda=\la^{(i)}$ for some row $i$) we have
\[c(\Lambda \setminus\la)=\la_i-(i-1)/\alpha=
\Lambda_i-1-(i-1)/\alpha.\]
Hence the moments of the random variable $c$ with respect to 
the transition and co-transition distributions are 
respectively
\begin{equation*}
\begin{split}
M_{r}(c)&= \sum_{\Lambda :\la \nearrow \Lambda} 
{c(\Lambda \setminus\la)}^r \,p_{\la}(\Lambda)=s_r(\la),\\
\tilde{M}_{r}(c)&= \sum_{\la :\la \nearrow \Lambda}
{c(\Lambda \setminus\la)}^r \,q_{\Lambda}(\la)=
\frac{1}{|\Lambda|}
\sum_{k=0}^r (-1)^{r-k} \binom{r}{k} \sigma_k(\Lambda).
\end{split}
\end{equation*}
Our results in Sections 8 and 9 thus amount to give an 
explicit evaluation of these moments.

If $\alpha =1$ the multiplicities are given by $\kappa(\la,\Lambda)=
H_{\la}/H_{\Lambda}$, with $H_{\la}$ the classical hook 
polynomial $\prod_{i,j \in \la}(\la_i+\la^{'}_j-i-j+1)$, 
and $\la'$ the partition conjugate to $\la$. 
Then $H_{\la}\ \mathrm{dim}\ 
\la=f_{\la}$, with $f_{\la}$ the dimension of the irreducible 
representation of the symmetric group associated to $\la$ (i.e. the number of 
standard tableaux of shape $\la$). Thus one recovers 
the classical distributions
\[p_{\la}(\Lambda)= \frac{1}{|\Lambda|}
\frac{f_{\Lambda}}{f_{\la}}\quad,\quad
q_{\Lambda}(\la)= 
\frac{f_{\la}}{f_{\Lambda}}\]
Our results are new even in this classical case.

In~\cite{Ke} it was shown that Young diagrams must be interpreted as a 
special case of ``pairs of interlacing sequences'', and that 
transition and co-transition distributions can be defined 
in this more general frame. It is likely that our explicit evaluations 
can be easily translated in this context.

\section{Rows and columns}

To conclude we give the proof of two conjectures 
stated in Sections 9 and 10 of~\cite{La5} (all conjectures 
of these sections will then be established). This result was not included 
in~\cite{LL}, though it is a direct consequence of relation 
(5.1) proved there. We shall need the following lemma.

\begin{prop}
With $s(n,k)$ the Stirling numbers of the first kind, we have
\[\frac{1}{x^k}=\sum_{n\ge k} s(n-1,k-1)\frac{1}{[x]_n}.\]
\end{prop}
\begin{proof}
Let $\Delta$ be the finite difference operator defined by 
$\Delta f(y)=f(y+1)-f(y)$. We write the Newton interpolation formula
\[f(y)=\sum_{n\ge 0}[y]_n \, \frac{\Delta^n f(0)}{n!},\]
for $f(y)=x/(x-y)=\sum_{r\ge 0} y^r/x^r$.
This yields immediatly
\[f(y)=\sum_{n\ge 0} [y]_n \, \frac{x}{[x]_{n+1}}.\]
Hence the result by identifying coefficients of $y^k$.
\end{proof}

\begin{theorem}
For any integer $p\ge 0$ we have
\begin{equation*}
\begin{split}
{\binom{\la}{(p)}}_{\alpha} = \frac{1}{(1/\alpha)_p} 
\sum_{0 <i+j \le p}
&\frac{1}{\alpha^i} \, s(p-1,i+j-1)
\left(\sum_{k=0}^{min(i,j)}\binom{|\la|-j}{i-k}\,F_{jk}(\la) \right)\\
{\binom{\la}{1^p}}_{\alpha} = \frac{1}{(\alpha)_p} 
\sum_{0 <i+j \le p} 
(-1)^j& \alpha^{i+j} s(p-1,i+j-1)
\left(\sum_{k=0}^{min(i,j)}\binom{|\la|-j}{i-k}\,F_{jk}(\la) \right).
\end{split}
\end{equation*}
\end{theorem}
\begin{proof}
By the property of duality (\cite{La6}, p. 320)
\[{\binom{\la}{(p)}}_{\alpha} ={\binom{\la'}{1^p}}_{1/\alpha}\]
and using $d_{k}(\la';1/\alpha)=(-\alpha)^k d_{k}(\la;\alpha)$, 
hence $F_{jk}(\la';1/\alpha)=(-\alpha)^j F_{jk}(\la;\alpha)$,
it is enough to prove the second formula. The generalized Chu-Vandermonde 
formula proved in~\cite{La5} (Theorem 12.1, p. 161) reads
\begin{equation*}
\frac{(y+1)_{\la}}{(y)_{\la}} =
\sum_{\mu} (-1)^{\mu} {\binom{\la}{\mu}}_{\alpha} 
\frac{(-1)_\mu}{(y)_\mu}=
\sum_{p\ge 0} {\binom{\la}{1^p}}_{\alpha} \, 
\frac{(\alpha)_p}{[\alpha y]_p},
\end{equation*}
since $(-1)_\mu=0$ if $\mu_{1} > 1$.
On the other hand relation (5.1) implies
\begin{equation*}
\frac{(y+1)_{\la}}{(y)_{\la}} =
\sum_{i,j\ge 0} (-1)^{j} \frac{1}{y^{i+j}}
\left(\sum_{k=0}^{min(i,j)}\binom{|\la|-j}{i-k}\,F_{jk}(\la) \right).
\end{equation*}
Using Lemma 11.1, we conclude by comparison.
\end{proof}

\section*{Acknowledgements}
It is a pleasure to thank Alain Lascoux for pointing out
the relevance of Lagrange interpolation.

\end{document}